\documentclass{amsart}
\usepackage{amsfonts, amsmath, latexsym, epsfig}
\usepackage[norelsize]{algorithm2e}
\usepackage{amssymb}
\usepackage{color}
\usepackage{epsf}
\usepackage{url}

\newcommand{\RR}{\ensuremath{\mathbb{R}}}

\newcommand{\ZZ}{\ensuremath{\mathbb{Z}}}

\newtheorem{proposition}{Proposition}
\newtheorem{theorem}{Theorem}

\def\QuotS#1#2{\leavevmode\kern-.0em\raise.2ex\hbox{$#1$}\kern-.1em/\kern-.1em\lower.25ex\hbox{$#2$}}

\DeclareMathOperator{\GL}{GL}

\begin{document}

\author{Michel Deza}
\address{Michel Deza, \'Ecole Normale Sup\'erieure, 75005 Paris}
\email{Michel.Deza@ens.fr}

\author{Mathieu Dutour Sikiri\'c}
\address{Mathieu Dutour Sikiri\'c, Rudjer Boskovi\'c Institute, Bijenicka 54, 10000 Zagreb, Croatia, Fax: +385-1-468-0245}
\email{mathieu.dutour@gmail.com}

\thanks{Second author gratefully acknowledges support from the Alexander von Humboldt foundation}

\title{The hypermetric cone on $8$ vertices and some generalizations}

\date{}

\maketitle

\begin{abstract}
The lists of facets -- $298,592$  in $86$ orbits -- and of extreme rays -- $242,695,427$ in $9,003$ orbits --
of the  hypermetric cone $HYP_8$ are computed. 
The first generalization considered is the hypermetric polytope $HYPP_n$ for which we give general algorithms and a description for $n\le 8$.
Then we shortly consider generalizations to simplices of volume higher than $1$, hypermetric on graphs and infinite dimensional hypermetrics.
\end{abstract}

\section{Introduction}
 
Metric, cut and hypermetric cones
 are among central objects of Discrete Mathematics.
For example, finite metrics and $l_1$-metrics can be studied by polyhedral cones or polytopes; see, for example, \cite{DL}.
The {\em hypermetric cone} $HYP_n$ is the set of all {\em hypermetrics on $n$ points}, i.e., the  functions $d:\{1,\dots,n\}^2 \rightarrow \RR$, such that
\begin{equation*}
H(b, d)=\sum_{1\leq i<j\leq n} b_i b_j d(i,j) \leq 0 \mbox{~for~all~} b\in \ZZ^n, \sum_i b_i=1.
\end{equation*}
If   $b$ is $\{0,\pm 1\}$-valued and has
$k+1$ ones,  the inequality is called 
{\em $(2k+1)$-gonal inequality}. In fact, the case of general $b$ can be seen as some such
$\{0,\pm 1\}$-valued $b$ on a multiset of  $\sum_{i=1}^n|b_i|$ points, in which different
points occur $|b_1|,\dots , |b_n|$ times.

The {\em metric cone} $MET_n$ is  the set of all {\em semimetrics on $n$ points}, i.e., those of above functions $d$,
which  satisfy all {\em triangle} (i.e., all $3$-gonal) {\em inequalities}.

For a set $S\subseteq \{1, \dots, n\}$ the {\em cut} (or {\em split}) {\em semimetric} $\delta_S$ is a vector 
 defined as
\begin{equation*}
\delta_S(x,y)=\left\{\begin{array}{rcl}
1 &\mbox{~~~if~~~}  \vert S \cap \{x,y\}\vert = 1\\
0 &\mbox{otherwise.}
\end{array}\right.
\end{equation*}
Clearly, $\delta_{\overline{S}}= \delta_S$, and
$\delta_S$ is also can be seen as the adjacency matrix of a {\em cut}  (into $S$ and $\overline{S}$) {\em subgraph} of $K_n$. 
The {\em cut cone} $CUT_n$ is the positive span of the  $2^{n-1}- 1$ non-zero cut semimetrics;
 the {\em cut polytope} $CUTP_n$ is  the convex hull of all
$2^{n-1}$ cut semimetrics.

We have the evident inclusions $CUT_n \subseteq HYP_n \subseteq MET_n$  with 
 $CUT_n = MET_n$ only for $3 \le n \le 4$; also, $CUT_n = HYP_n$ only for $3 \le n \le 6$ (\cite{DL}).
 So, the first proper $HYP_n$ is $HYP_7$; it was described in \cite{Hyp7}.
While $CUT_n$ is important in Analysis and Combinatorics, the cone $HYP_n$ is related
to Voronoi-Delaunay theory in Geometry of Numbers.
 
 Given a polyhedral cone  $C$, let $e_C$ denote the number of extreme rays of $C$, while
 $f_C$ denote the number of facets of  $C$.
The Table \ref{tab1} gives $e_C$ and $f_C$ for $C=CUT_n,HYP_n,MET_n$ with $n=3,4,5,6,7,8$.
The enumeration of orbits of facets of $CUT_n$  for $n\leq 7$ was done in \cite{OrbitFacetCutPolytope5,OrbitFacetCutPolytope6,OrbitFacetCutPolytope7} for $n=5$, $6$ and $7$, respectively.
For $CUT_8$ and $CUTP_8$, sets of facets  were found in \cite{CR};
completeness of these sets was shown in \cite{DD15}.
The enumeration of orbits of extreme rays of $MET_n$ for $n\leq 8$ was done in \cite{OrbitRaysMetricCone7,OrbitVerticesMetricPolytope7,OrbitPolytopeMetricPolytope8}.

\begin{table}
\begin{center}
\caption{The number of extreme rays  and facets in cones $HYP_n$, $CUT_n$ and $MET_n$ for $3\le n\le 8$; the numbers of orbits under the symmetric group $Sym(n)$ are given in parentheses. Also, the number of vertices  and facets in the hypermetric polytope $HYPP_n$, see Section 4, is given for $3\le n\le 8$, with number of orbits under  $Sym(n)$ and $2^{n-1}$ switchings}

\label{tab1}
\begin{tabular}{||c||c|c|c|c|c|c||}
\hline
\hline
$C$              & $n=3$  & $n=4$   & $n=5$   & $n=6$   & $n=7$             & $n=8$\\
\hline
\hline
$CUT_n, e$       & $3(1)$ & $7(2)$  & $15(2)$ & $31(3)$ & $63(3)$           & $127(4)$\\
$CUT_n,f$        & $3(1)$ & $12(1)$ & $40(2)$ & $210(4)$& $38,780(36)$      & $49,604,520(2,169)$\\
\hline
${\bf HYP_n,e}$  & $3(1)$ & $7(2)$  & $15(2)$ & $31(3)$ & $\bf{37,170(29)}$ & $\bf{242,695,427(9,003)}$\\
${\bf HYP_n,f}$  & $3(1)$ & $12(1)$ & $40(2)$ & $210(4)$& $\bf{3,773(14)}$  & $\bf{298,592(86)}$\\
\hline
$MET_n,e$        & $3(1)$ & $7(2)$  & $25(3)$ & $296(7)$& $55,226(46)$      & $119,269,588(3,918)$\\
$MET_n,f$        & $3(1)$ & $12(1)$ & $30(1)$ & $60(1)$ & $105(1)$          & $168(1)$\\
\hline
\hline
${\bf HYPP_n,v}$ & $4(1)$ & $8(1)$  & $16(1)$ & $32(1)$ & $\bf{113,152(6)}$ &$\bf{1,388,383,872(581)}$\\
${\bf HYPP_n,f}$ & $4(1)$ & $16(1)$ & $56(2)$ & $68(3)$ & $\bf{10,396(7)}$  &$\bf{1,374,560(22)}$\\
\hline
\hline
\end{tabular}
\end{center}
\end{table}

In Section \ref{Facet_HYP8} the facets of the hypermetric cone $HYP_8$ are determined with the help of the connection with geometry of numbers and the list of simplices of dimension $7$.
In Section \ref{ExtremeRay_HYP8} we list the extreme rays of $HYP_8$ by using the list of extreme Delaunay polytopes in dimension $7$.

In Section \ref{HypermetricPolytope} we define the hypermetric polytope and give algorithms for computing with it. This is then used to compute the vertices and facets of the hypermetric polytope $HYPP_7$ and $HYPP_8$.
In Section \ref{Structure_AllSimplices} we give the list of facets of the cone that come from simplices of volume greater than $1$. Those are direct analogues of the hypermetric cone.
In Section \ref{Hypermetric_OnGraphs} we define an analogue of hypermetric inequality for the cut-polytope of a graph $G$. This directly generalizes the corresponding definitions of the metric polytope of a graph and allows to find new facet inequality in some cases.
In Section \ref{Infinite_Hypermetric} the notion of infinite hypermetric is briefly considered.

\section{Computation of facets of $HYP_8$}\label{Facet_HYP8}

While $HYP_n$ is defined by an infinity of inequalities, it is proved 
in \cite{DGL92,DGL95,DL}
that this cone is, in fact, polyhedral for all $n$. The proof relies on a connection with Geometry of Numbers that we now explain.

Given a quadratic form $q$, one can define the induced Delaunay tessellation with point set $\ZZ^n$ (\cite{VoronoiII,bookschurmann,EquivariantLtypeDSV,ComplexityVoronoiDSV}).
It is well known, that in this context there are only a finite number of possible tessellations, up to the action of the group $\GL_n(\ZZ)$. 

For a generic quadratic form, the tessellation is formed by simplices only; but, importantly, when it is not, this induces linear conditions on the coefficients of the 
 form.
There are a finite number of simplices, up to $\GL_n(\ZZ)$ action, and they have been classified in \cite{InhomogeneousPerfect} for
 $n=7$, extending previous classification for
 $n\le 6$ (\cite{Ctype_original,Ba,BaRy}).
There are $11$ orbits of such simplices and their volume is at most $5$.

Given a simplex $S$, one can consider all Delaunay polytopes containing it. We consider the set $Bar_S$ of quadratic forms, for which $S$ is contained in a Delaunay polytope of the Delaunay tessellation. 
It is a polyhedral cone,  called  a {\em Baranovskii cone} in \cite{bookschurmann}.

For a quadratic form inside $Bar_S$, the Delaunay tessellation contains $S$ as a simplex. For a quadratic form on a facet of $Bar_S$, the simplex $S$ is a part of a  {\em repartitioning polytope}, i.e., a Delaunay polytope with $n+2$ vertices.

If the simplex $S$ has volume $1$, then it is equivalent to the simplex  formed by the vertices $v_1=0$, $v_2=e_1$, \dots, $v_{n+1}=e_n$.
The quadratic form $q$ is described uniquely by the distance function $d(i,j)=q(v_i - v_j)$ on the vertices $v_i$.
For a given positive definite quadratic form $q$, 
 denote by $c(q)$ the center of the sphere, circumscribing $S$, and by $r(q)$ the radius of this sphere.
Since $S$ is of volume $1$, a given point $v\in \ZZ^n$ can be uniquely expressed in barycentric coordinates as $v=\sum_i b_i v_i$ with $1=\sum_i b_i$ and $b_i\in \ZZ$.
In \cite{DL}, the following formula is proved:
\begin{equation*}
H(b, d)=\Vert v - c(q) \Vert^2 - r(q)^2.
\end{equation*}
So, the distance function $d$ corresponds to a quadratic form $q$ having $S$ a part of the Delaunay tessellation if and only if it belongs to $HYP_{n+1}$.
Therefore, in order to classify the facets of $HYP_8$, we need to classify the repartitioning polytopes 
of dimension 
 $7$.

A repartitioning polytope $R$ is defined by its $n+2$ vertices $w_1, \dots, w_{n+2}$. There exists a unique linear relation among the vertices. It is expressed as
\begin{equation*}
\sum_{i=1}^{n+2} \alpha_i w_i = 0 \mbox{~with~}\sum_{i=1}^{n+2} w_i = 0.
\end{equation*}
It turns out, that $R$ admits exactly two triangulations into simplices. Define
\begin{equation*}
S_+=\{1\leq i\leq n+2\mbox{~s.t.~} w_i > 0\} \mbox{~and~}S_-=\{1\leq i\leq n+2\mbox{~s.t.~} w_i < 0\}.
\end{equation*}
The first triangulation is defined by taking all simplices with vertices $\{1,\dots, n+2\} - \{i\}$ for $i\in S_+$ and the second one similarly from $S_-$. Going from one triangulation to another is called a {\em bistellar flips} \cite{BistellarFlip}.

A $L$-type domain is the specification of the full tessellation by Delaunay polytopes of $\ZZ^n$.
The $L$-type form a tessellation of the cone of positive definite quadratic forms.
When one moves from one $L$-type to another adjacent one then the operation that is done is a change
from one triangulation $T$ to another $T'$. This change can be described in the following way:
\begin{enumerate}
\item Some of the simplices of the triangulation $T$ are grouped to form repartitioning polytopes
\item Each repartitioning polytope has two triangulations. We go from one to another.
\item The triangulation $T'$ is obtained from all the new triangulations and the simplices that did not
belong to any repartitioning polytope.
\end{enumerate}
As a consequence, all simplices occurring in a repartitioning polytope, are also lattice polytopes.

Now we explain our strategy for enumerating the repartitioning polytopes, which is rather similar to the one of \cite{SmoothnessRegularity}:
\begin{theorem}
There are exactly $67$ types of repartitioning polytopes for $n=7$.
\end{theorem}
\proof Our classification is based on a computer assisted case distinction.
The full list can be obtained from the web-page~\cite{WebPageRepartitioningDim7}.
Here, we briefly describe the necessary ingredients and our computational steps.

Let us write the vertices of a repartitioning polytope $R$ as $v_1, \dots, v_9$.
Without loss of generality, assume that $v_1, \dots, v_8$ are linearly independent.

For any $1\leq i\leq n$, we define $S_i=\{v_1, \dots, v_9\} - \{v_i\}$ and denote
by $vol(S_i)$ its volume multiplied by $n!$.
If the points of $S_i$ are linearly dependent, then $vol(S_i)=0$; otherwise,
$S_i$ is a realizable simplex and so, by the classification of \cite{InhomogeneousPerfect},
of volume at most $5$.
In addition, one can assume that $vol(S_9)$ is maximal among all $vol(S_i)$ and
denote it by $i(R)\leq 5$.
By using exterior algebra product, one can find vectors $w_i$ such that
\begin{equation*}
vol(S_i) = \left\vert \langle w_i, v_9\rangle \right\vert.
\end{equation*}
Since either $vol(S_i)=0$, or $vol(S_i)\leq vol(S_9)=i(R)$, we have the inequalities:
\begin{equation*}
-i(R) \leq \langle w_i, v_9\rangle \leq i(R).
\end{equation*}

Geometrically, these conditions define a polytope, for which we are
searching its integral points.
 One can use, for instance, the program {\tt zsolve} 
from \cite{4ti2} to enumerate those integral points.
The polytope thus defined is a little more complicate than a parallelepiped;
but, as $i(R)$ increases, the condition  $vol(S_9)=i(R)$
impose some strong restriction on the solution set, which, so,
does not increases too much.

Hence, we obtain a finite list of possible candidates for the
repartitioning polytopes.
Then we need to check whether there
exist an adequate quadratic form realizing it. This is done by an
adaptation of Algorithm~1 of~\cite{MinkowskianSublattices}. That is,
we use the integral symmetries of the repartitioning polytopes and
iterate until an adequate quadratic form is found. \qed

After obtaining the list of $67$ repartitioning polytopes, we look
at all the simplices of volume $1$ in it and  at the corresponding
barycentric coordinates of the remaining vertex. Therefore, we get:

\begin{theorem}
The hypermetric cone $HYP_8$ has $298,592$ facets in $86$ orbits.
\end{theorem}
Note that  \cite{DutourAdj} gave $86$ as lower bound on the
number of orbits of facets, which is therefore an exact bound.

The $86$ orbits $O_i$ of facets are presented in Tables \ref{tab2_1}, \ref{tab2_2}.
The first representatives of each (of $22$) switching equivalence classes of facets are boldfaced there.
The orbits of simplicial facets are marked by $^{*}$. About $92\%$ of the total number of facets ($60$ orbits) are simplicial; they polish the cone.
On the other hand, each triangle facet contains about $18.8$ millions of extreme rays.

For facets  in $HYP_n$,   gcd of orbit sizes
is  $56$ if $n=8$ and $3,12,10,30,7$ if $n=3,4,5,6,7$, respectively.

\begin{table}
\begin{center}
\caption{The orbits of facets of the cone $HYP_8$: part 1}
\label{tab2_1}
\begin{tabular}{||c|c|c|c|c||}
\hline
\hline
$F_{i,j}$ & Representative & $\frac{\left\vert F_{i,j}\right\vert}{56}$ & $CUT_8$-rank & Inc.($[0,1]$, $2_{21}$, $3_{21}$, $ER_7$)\\
\hline
\hline
$F_{1,1}$ & $\bf{(0,0,0,0,0,-1,1,1)}$ & $3$ & $27$ & $(95,329734,737128,17725428)$\\
\hline
\hline
$F_{2,1}$ & $\bf{(0,0,0,-1,-1,1,1,1)}$ & $10$ & $27$ & $(79,93978,176058,3780630)$\\
\hline
\hline
$F_{3,1}$ & $\bf{(0,0,-1,-1,-1,1,1,2)}$ & $30$ & $27$ & $(59,10460,13052,209644)$\\
$F_{3,2}$ & $(0,0,-1,1,1,1,1,-2)$ & $15$ & $27$ & $(59,10460,13052,209644)$\\
\hline
\hline
$F_{4,1}$ & $\bf{(0,-1,-1,-1,1,1,1,1)}$ & $5$ & $27$ & $(69,36816,60480,1207584)$\\
\hline
$F_{5,1}$ & $\bf{(0,-1,-1,-1,-1,1,1,3)}$ & $15$ & $27$ & $(41,400,240,620)$\\
$F_{5,2}$ & $(0,-1,1,1,1,1,1,-3)$ & $6$ & $27$ & $(41,400,240,620)$\\
\hline
$F_{6,1}$ & $\bf{(0,-1,-1,1,1,1,-2,2)}$ & $60$ & $27$ & $(51,3567,3288,46176)$\\
$F_{6,2}$ & $(0,-1,-1,-1,-1,1,2,2)$ & $15$ & $27$ & $(51,3650,4680,64400)$\\
$F_{6,3}$ & $(0,1,1,1,1,1,-2,-2)$ & $3$ & $27$ & $(51,3650,4680,64400)$\\
\hline
$F_{7,1}$ & $\bf{(0,-1,-1,-1,1,-2,2,3)}$ & $120$ & $26$ & $(39,311,220,1479)$\\
$F_{7,2}$ & $(0,-1,-1,1,1,2,2,-3)$ & $90$ & $26$ & $(39,325,172,1461)$\\
$F_{7,3}$ & $(0,-1,1,1,1,-2,-2,3)$ & $60$ & $26$ & $(39,325,172,1461)$\\
$F_{7,4}$ & $(0,1,1,1,1,-2,2,-3)$ & $30$ & $26$ & $(39,311,220,1479)$\\
\hline
\hline
$F_{8,1}$ & $\bf{(-1,-1,-1,-1,1,1,1,2)}$ & $5$ & $27$ & $(55,6840,8526,141642)$\\
$F_{8,2}$ & $(-1,-1,1,1,1,1,1,-2)$ & $3$ & $27$ & $(55,6840,8526,141642)$\\
\hline
$F_{9,1}$ & $\bf{(-1,-1,-1,-1,-1,1,1,4)}$ & $3$ & $27$ & $(27,0,0,0)^*$\\
$F_{9,2}$ & $(-1,1,1,1,1,1,1,-4)$ & $1$ & $27$ & $(27,0,0,0)^*$\\
\hline
$F_{10,1}$ & $\bf{(-1,-1,-1,1,1,1,-2,3)}$ & $20$ & $27$ & $(41,645,282,3021)$\\
$F_{10,2}$ & $(-1,-1,1,1,1,1,2,-3)$ & $15$ & $27$ & $(41,645,282,3021)$\\
$F_{10,3}$ & $(-1,-1,-1,-1,-1,1,2,3)$ & $6$ & $27$ & $(41,495,828,5094)$\\
$F_{10,4}$ & $(1,1,1,1,1,1,-2,-3)$ & $1$ & $27$ & $(41,495,828,5094)$\\
\hline
$F_{11,1}$ & $\bf{(-1,-1,-1,1,1,-2,2,2)}$ & $30$ & $27$ & $(45,1464,1390,18310)$\\
$F_{11,2}$ & $(-1,1,1,1,1,-2,-2,2)$ & $15$ & $27$ & $(45,1464,1390,18310)$\\
$F_{11,3}$ & $(-1,-1,-1,-1,-1,2,2,2)$ & $1$ & $27$ & $(45,2070,1458,34956)$\\
\hline
$F_{12,1}$ & $\bf{(-1,-1,1,1,-2,-2,2,3)}$ & $90$ & $27$ & $(37,293,166,1638)$\\
$F_{12,2}$ & $(-1,1,1,1,-2,2,2,-3)$ & $60$ & $27$ & $(37,293,166,1638)$\\
$F_{12,3}$ & $(-1,-1,-1,1,2,2,2,-3)$ & $20$ & $27$ & $(37,306,279,2616)$\\
$F_{12,4}$ & $(-1,-1,-1,-1,-2,2,2,3)$ & $15$ & $27$ & $(37,300,285,2500)$\\
$F_{12,5}$ & $(1,1,1,1,-2,-2,-2,3)$ & $5$ & $27$ & $(37,306,279,2616)$\\
\hline
$F_{13,1}$ & $\bf{(-1,-1,-1,-1,1,-2,2,4)}$ & $30$ & $26$ & $(31,35,31,31)$\\
$F_{13,2}$ & $(-1,-1,1,1,1,-2,-2,4)$ & $30$ & $26$ & $(31,45,3,24)$\\
$F_{13,3}$ & $(-1,-1,1,1,1,2,2,-4)$ & $30$ & $26$ & $(31,45,3,24)$\\
$F_{13,4}$ & $(1,1,1,1,1,-2,2,-4)$ & $6$ & $26$ & $(31,35,31,31)$\\
\hline
$F_{14,1}$ & $\bf{(-1,-1,-1,1,1,2,-3,3)}$ & $60$ & $26$ & $(35,142,46,268)$\\
$F_{14,2}$ & $(-1,1,1,1,1,-2,-3,3)$ & $30$ & $26$ & $(35,142,46,268)$\\
$F_{14,3}$ & $(-1,-1,-1,-1,1,-2,3,3)$ & $15$ & $26$ & $(35,110,142,404)$\\
$F_{14,4}$ & $(1,1,1,1,1,2,-3,-3)$ & $3$ & $26$ & $(35,110,142,404)$\\
\hline
$F_{15,1}$ & $\bf{(-1,-1,-1,-1,1,-3,3,4)}$ & $30$ & $26$ & $(26,0,0,1)^*$\\
$F_{15,2}$ & $(-1,-1,-1,1,1,3,3,-4)$ & $30$ & $26$ & $(26,0,0,1)^*$\\
$F_{15,3}$ & $(-1,1,1,1,1,-3,-3,4)$ & $15$ & $26$ & $(26,0,0,1)^*$\\
$F_{15,4}$ & $(1,1,1,1,1,-3,3,-4)$ & $6$ & $26$ & $(26,0,0,1)^*$\\
\hline\hline
\end{tabular}
\end{center}
\end{table}

\begin{table}
\begin{center}
\caption{The orbits of facets of the cone $HYP_8$: part 2}
\label{tab2_2}
\begin{tabular}{||c|c|c|c|c||}
\hline
\hline
$F_{i,j}$ & Representative & $\frac{\left\vert F_{i,j}\right\vert}{56}$ & $CUT_8$-rank & Inc.($[0,1]$, $2_{21}$, $3_{21}$, $ER_7$)\\
\hline
\hline
$F_{16,1}$ & $\bf{(-1,-1,1,-2,2,2,-3,3)}$ & $180$ & $26$ & $(32,63,36,177)$\\
$F_{16,2}$ & $(1,1,1,-2,-2,2,-3,3)$ & $60$ & $26$ & $(32,63,36,177)$\\
$F_{16,3}$ & $(-1,1,1,2,2,2,-3,-3)$ & $30$ & $26$ & $(32,82,38,272)$\\
$F_{16,4}$ & $(-1,1,1,-2,-2,-2,3,3)$ & $30$ & $26$ & $(32,82,38,272)$\\
$F_{16,5}$ & $(-1,-1,-1,-2,-2,2,3,3)$ & $30$ & $26$ & $(32,94,26,266)$\\
\hline
$F_{17,1}$ & $\bf{(-1,-1,1,1,-2,2,-3,4)}$ & $180$ & $25$ & $(30,41,6,51)$\\
$F_{17,2}$ & $(-1,1,1,1,-2,2,3,-4)$ & $120$ & $25$ & $(30,41,6,51)$\\
$F_{17,3}$ & $(-1,-1,-1,1,2,2,3,-4)$ & $60$ & $25$ & $(30,40,22,74)$\\
$F_{17,4}$ & $(-1,-1,-1,1,-2,-2,3,4)$ & $60$ & $25$ & $(30,30,32,74)$\\
$F_{17,5}$ & $(1,1,1,1,-2,-2,-3,4)$ & $15$ & $25$ & $(30,40,22,74)$\\
$F_{17,6}$ & $(1,1,1,1,2,2,-3,-4)$ & $15$ & $25$ & $(30,30,32,74)$\\
$F_{17,7}$ & $(-1,-1,-1,-1,2,2,-3,4)$ & $15$ & $25$ & $(30,40,22,72)$\\
\hline
$F_{18,1}$ & $\bf{(-1,-1,1,-2,2,3,3,-4)}$ & $180$ & $25$ & $(27,13,9,20)$\\
$F_{18,2}$ & $(-1,1,1,2,2,-3,3,-4)$ & $180$ & $25$ & $(27,15,3,15)$\\
$F_{18,3}$ & $(-1,1,1,-2,-2,-3,3,4)$ & $180$ & $25$ & $(27,15,3,15)$\\
$F_{18,4}$ & $(-1,-1,-1,-2,2,-3,3,4)$ & $120$ & $25$ & $(27,13,9,18)$\\
$F_{18,5}$ & $(-1,-1,1,2,2,-3,-3,4)$ & $90$ & $25$ & $(27,21,1,22)$\\
$F_{18,6}$ & $(1,1,1,-2,2,-3,-3,4)$ & $60$ & $25$ & $(27,13,9,20)$\\
$F_{18,7}$ & $(1,1,1,-2,-2,3,3,-4)$ & $30$ & $25$ & $(27,21,1,22)$\\
\hline
$F_{19,1}$ & $\bf{(-1,-1,-1,1,-2,-2,2,5)}$ & $60$ & $24$ & $(24,2,1,0)^*$\\
$F_{19,2}$ & $(-1,-1,1,1,2,2,2,-5)$ & $30$ & $24$ & $(24,3,0,0)^*$\\
$F_{19,3}$ & $(-1,1,1,1,-2,-2,-2,5)$ & $20$ & $24$ & $(24,3,0,0)^*$\\
$F_{19,4}$ & $(1,1,1,1,-2,2,2,-5)$ & $15$ & $24$ & $(24,2,1,0)^*$\\
\hline
$F_{20,1}$ & $\bf{(-1,-1,-1,1,-2,-3,3,5)}$ & $120$ & $24$ & $(24,1,1,1)^*$\\
$F_{20,2}$ & $(-1,-1,1,1,2,-3,-3,5)$ & $90$ & $24$ & $(24,2,0,1)^*$\\
$F_{20,3}$ & $(-1,-1,-1,1,2,3,3,-5)$ & $60$ & $24$ & $(24,2,0,1)^*$\\
$F_{20,4}$ & $(-1,1,1,1,-2,3,3,-5)$ & $60$ & $24$ & $(24,2,0,1)^*$\\
$F_{20,5}$ & $(1,1,1,1,2,-3,3,-5)$ & $30$ & $24$ & $(24,1,1,1)^*$\\
$F_{20,6}$ & $(1,1,1,1,-2,-3,-3,5)$ & $15$ & $24$ & $(24,2,0,1)^*$\\
\hline
$F_{21,1}$ & $\bf{(-1,1,-2,2,2,-3,-3,5)}$ & $180$ & $23$ & $(23,2,1,1)^*$\\
$F_{21,2}$ & $(-1,-1,-2,-2,2,-3,3,5)$ & $180$ & $23$ & $(23,3,0,1)^*$\\
$F_{21,3}$ & $(-1,1,2,2,2,-3,3,-5)$ & $120$ & $23$ & $(23,2,1,1)^*$\\
$F_{21,4}$ & $(1,1,-2,-2,2,3,3,-5)$ & $90$ & $23$ & $(23,2,1,1)^*$\\
$F_{21,5}$ & $(-1,-1,-2,2,2,3,3,-5)$ & $90$ & $23$ & $(23,3,0,1)^*$\\
$F_{21,6}$ & $(1,1,-2,-2,-2,-3,3,5)$ & $60$ & $23$ & $(23,2,1,1)^*$\\
\hline
$F_{22,1}$ & $\bf{(-1,-1,1,-2,2,3,4,-5)}$ & $360$ & $23$ & $(23,2,1,1)^*$\\
$F_{22,2}$ & $(-1,-1,1,-2,-2,-3,4,5)$ & $180$ & $23$ & $(23,2,1,1)^*$\\
$F_{22,3}$ & $(-1,1,1,-2,-2,3,-4,5)$ & $180$ & $23$ & $(23,3,0,1)^*$\\
$F_{22,4}$ & $(-1,1,1,2,2,-3,4,-5)$ & $180$ & $23$ & $(23,3,0,1)^*$\\
$F_{22,5}$ & $(-1,-1,1,2,2,-3,-4,5)$ & $180$ & $23$ & $(23,3,0,1)^*$\\
$F_{22,6}$ & $(1,1,1,-2,2,-3,-4,5)$ & $120$ & $23$ & $(23,2,1,1)^*$\\
$F_{22,7}$ & $(-1,-1,-1,-2,2,3,-4,5)$ & $120$ & $23$ & $(23,2,1,1)^*$\\
$F_{22,8}$ & $(1,1,1,2,2,3,-4,-5)$ & $60$ & $23$ & $(23,2,1,1)^*$\\
$F_{22,9}$ & $(1,1,1,-2,-2,3,4,-5)$ & $60$ & $23$ & $(23,3,0,1)^*$\\
\hline\hline
\end{tabular}
\end{center}
\end{table}

\section{Extreme rays of $HYP_8$}\label{ExtremeRay_HYP8}
A $n$-dimensional Delaunay polytope $D$ is called {\em extreme} if up to scalar multiple there is a unique quadratic
form $q$ having $D$ as Delaunay polytopes.
For a Delaunay polytope $D$ an {\em integral affine generating set} $S_{aff}$ is a set $v_1, \dots, v_m$ of vertices
of $D$ such that for each vertex $v$ of $D$ there exist $\lambda_i\in \ZZ$ such that
\begin{equation*}
v = \sum_{i=1}^m \lambda_i v_i \mbox{~and~} 1 = \sum_{i=1}^m \lambda_i.
\end{equation*}
In the special case $m=n+1$ the set $S_{aff}$ is called an {\em affine basis}.
Given an extreme Delaunay polytope $D$ of associated quadratic form $q$ and an integral affine generating set
$S_{aff}=\left\{v_1, \dots, v_m\right\}$ the distance function $d(i,j) = q(v_i - v_j)$ defines an extreme ray
of $HYP_m$.
It is proved in \cite{DL} that all extreme rays of $HYP_n$ are obtained in this way.
Thus in order to classify the extreme rays of $HYP_n$ we need the list of extreme Delaunay polytopes of
dimension at most $n-1$ and their integral affine generating sets.

The perfect Delaunay polytopes of dimension at most $6$ are determined in \cite{Hyp7}.
The ones of dimension $7$ are determined in \cite{InhomogeneousPerfect}.
In summary, we obtain the following:
\begin{enumerate}
\item The $1$-dimensional polytope $[0,1]$. It has $4$ orbits of $8$-points integral affine generating sets and this yield $127$ extreme rays of $HYP_8$.
\item The $6$-dimensional {\em Schlafli polytope} $2_{21}$. It has $195$ orbits of $8$ points integral affine generating sets. This gives a total of $231,596$ extreme rays in $HYP_8$.
\item The $7$-dimensional {\em Gosset polytope} $3_{21}$. It has $374$ orbits of affine basis and gives a total of $7,126,560$ extreme rays of $HYP_8$.
\item The $7$-dimensional {\em Erdahl-Rybnikov polytope} $ER_7$.
It has $8,430$ orbits of affine basis and gives a total of $235,337,144$ extreme rays of $HYP_8$.

\end{enumerate}

For $P=3_{21}$ or $ER_{7}$, it suffices to enumerate the
orbits of $8$ vertices in $P$. We keep the ones that determine a simplex of volume $1$ and thus
are integral affinely generating. We found $374$ and $8,430$ orbits.

For $P=2_{21}$, the integral affine generating sets can have $7$ points with one repeated or $8$ points.
In the $7$-points case we enumerate the affine basis of $2_{21}$ and consider all ways to duplicate points.
In The $8$-points case we enumerate all orbits of $8$-points and check the ones that integrally affine
generates $2_{21}$. This gives $195$ orbits.

For the interval $[0,1]$, we simply have to look at the number $2^{n-1}-1$ of non-zero cut semimetrics
 on $n$ vertices. 
 
So, altogether $HYP_8$ has $242,695,427$ vertices in $9,003$ orbits.

\section{Hypermetric polytope}\label{HypermetricPolytope}

For both, the metric cone and the cut cone, there are polytope analogs.
The {\em cut polytope} $CUTP_n$ is defined as the convex hull of all
$2^{n-1}$ cut semimetrics.
The {\em metric polytope} $METP_n$ is defined by the same  $3{n\choose 3}$
triangle inequalities, as $MET_n$,
and ${n\choose 3}$ additional {\em perimeter inequalities}
\begin{equation*}
d_{ij} + d_{jk} + d_{ki} \leq 2 \mbox{~for~}1\leq i<j<k\leq n.
\end{equation*}

The polytopes $CUTP_n$ and $METP_n$ are invariant under the
following {\em switching} operation $U_S$ on semimetrics:
\begin{equation*}
U_S(d)=\left\{\begin{array}{rcl}
\{1, \dots, n\}^2 &\rightarrow& \RR\\
(i,j)             &\mapsto    &
         \left\{\begin{array}{rl}
         1 - d(i,j) &\mbox{~if~}\vert S\cap \{i,j\}\vert = 1\\
         d(i,j)    &\mbox{~otherwise.}
         \end{array}\right.
\end{array}\right.
\end{equation*}

By analogy with above, we proceed in the following way for hypermetrics.
Given a vector $b=(b_1, \dots, b_n)\in \ZZ^n$ with $\sum_{i=1}^n b_i = 2s + 1$
and $s\in \ZZ$, we define the {\em hypermetric polytope} $HYPP_n$ by the inequalities
\begin{equation*}
\sum_{1\leq i<j\leq n} b_i b_j d(i,j) \leq s(s+1).
\end{equation*}
One obtains $METP_n$ using only  $b$ of the form  $(1,1,-1, 0^{n-3})$
and $(1,1,1,0^{n-3})$.

\begin{theorem}\label{CheckTheoremBelonging}
Given a distance function $d$, there is an algorithm for testing if

(i) $d\in HYP_n$

(ii) $d\in HYPP_n$
\end{theorem}
\proof To check if $d\in HYP_n$ is to check if for all $b\in \ZZ^{n}$ with $1=\sum_i b_i$ we have
\begin{equation*}
\sum_{1\leq i<j\leq n} b_i b_j d_{ij}\leq 0
\end{equation*}
while to check if $d\in HYPP_n$ is to check if for all $b\in \ZZ^n$ with $\sum_i b_i$ odd we have
\begin{equation*}
\sum_{1\leq i<j\leq n} b_i b_j d_{ij}\leq s(s+1) \mbox{~with~} 2s+1 = \sum_i b_i.
\end{equation*}
Both questions can be reframed in terms of quadratic functions, i.e. functions that are sum of a quadratic form, linear form, and constant term. 
Given a quadratic function $f$, we need to check if there exist a point $x\in \ZZ^n$ such that $f(x)<0$.
By standard linear algebra rewriting the question becomes to check if for a positive definite quadratic form $A$, vector $c\in \RR^n$ and distance $d$ there exist a vector $x\in \ZZ^n$ such that $A[x - c] < d$. This is a Closest Vector Problem and there are algorithms for solving such questions \cite{FinckePohst}. \qed

The skeletons of $HYPP_n, HYP_n$ contain a clique consisting of all cuts or all non-zero cuts. We expect that any vertex is adjacent to a cut vertex (it holds for $n\le 8$); if true, it will imply that each of  above skeletons have diameter $3$.
The {\em ridge graphs} (i.e., skeletons of the duals) of 
the triangle/perimeter facets of $METP_n, MET_n$ with $n\ge 4$ have diameter $2$
(\cite{RidgeGraphMetric}).
We expect that any facet of  $HYPP_n$, $HYP_n$ is adjacent to 
a triangle/perimeter facet (it holds for $n\le 7$); if true, it will imply that the ridge graphs of 
$HYPP_n, HYP_n$ have diameter $4$.

\begin{theorem}
The entries of Table \ref{tab1} for the hypermetric polytopes $HYPP_7$ and $HYPP_8$ are valid.
\end{theorem}
\proof The hypermetric cone $HYP_8$ has $9,003$ orbits of extreme rays and $86$ orbits
of facets under $Sym(8)$. We consider the group of order $2^7 8!$, denote it by $ARes(K_8)$,
generated by $Sym(8)$ and the switchings. 
Under the group $Ares(K_8)$ the facets of $HYP_8$ generate $22$ orbits of facets of
$HYPP_8$. The extreme rays of $HYP_8$ are of the form $\lambda v$ with $v$ a generator.
For each extreme ray, we choose $\lambda$ to be the maximal value which defines
a vertex of the hypermetric polytope $HYPP_8$ by using Theorem \ref{CheckTheoremBelonging}.
After elimination of isomorphic pairs, this gives $581$ orbits of vertices of $HYPP_8$.

However, there could be more facets of $HYPP_8$: In principle it could happen that 
a vector $(b_1, \dots, b_n)$ is not incident to any cut and yet defines a facet of
$HYPP_8$.
There could be more vertices as well: If a vertex of $HYPP_8$ is not adjacent to any
cut, then it does not appear from the list of extreme rays of $HYP_8$.
For each of the $581$ orbits of vertices we compute the adjacent vertices by using
the list of $22$ orbits that we have. All the vertices found belong to the $581$ orbits,
which proves that both lists are complete.

The same method applies as well to $HYPP_7$. \qed

The hypermetric polytope $HYPP_8$ has $581$ orbits of vertices, which are in details:
\begin{itemize}
\item $1$ orbit $V_C$  of $128$ cuts. 
The stabilizer of a cut is isomorphic to
Sym(8); the number of classes is $4$. $V_C$ forms a clique; the cone of facets incident
to a cut is exactly the hypermetric cone $HYP_8$. Cuts are the only vertices of $HYPP_n$, having all  coordinates integral.

\item $24$ orbits corresponding to the Delaunay polytopes $2_{21}$ and $3_{21}$. 
Details on those orbits are given in Table \ref{HYP8polytope_23_21_vertices}.
The denominator of the coordinates is $3$ for all vertices.
\item $556$ orbits corresponding the extreme Delaunay polytope $ER_7$. The denominator
of the coordinates is $12$ for all vertices. The number of incident inequalities and adjacent
vertices is $28$ for each of them.
\end{itemize}

\begin{table}
\begin{center}
\caption{Orbits of facets of the hypermetric polytope $HYPP_8$}
\label{HYP8polytope_facet}
\begin{tabular}{||c|c|c|c|c||}
\hline
\hline
$F_{i}$ & Representative & $\frac{\left\vert F_i\right\vert}{
32}$ & \#classes & Inc.($[0,1]$, $\{2_{21}, 3_{21}\}$, $ER_7$)\\
\hline
$F_{1}$ & $(0,0,0,0,0,1,1,1)$ & $7$ & $2$ & $(96,1598784,80836608)$\\
$F_{2}$ & $(0,0,0,1,1,1,1,1)$ & $28$ & $3$ & $(80,383040,14300640)$\\
$F_{3}$ & $(0,1,1,1,1,1,1,1)$ & $16$ & $4$ & $(70,131712,3975552)$\\
$F_{4}$ & $(0,0,1,1,1,1,1,2)$ & $168$ & $6$ & $(60,32160,590960)$\\
$F_{5}$ & $(0,1,1,1,1,1,2,2)$ & $336$ & $9$ & $(52,9600,122160)$\\
$F_{6}$ & $(1,1,1,1,1,1,1,2)$ & $32$ & $8$ & $(56,19656,370272)$\\
$F_{7}$ & $(0,1,1,1,1,1,1,3)$ & $112$ & $7$ & $(42,840,1120)$\\
$F_{8}$ & $(1,1,1,1,1,2,2,2)$ & $224$ & $12$ & $(46,3528,39906)$\\
$F_{9}$ & $(0,1,1,1,1,2,2,3)$ & $1,680$ & $15$ & $(40,656,2686)$\\
$F_{10}$ & $(1,1,1,1,1,1,2,3)$ & $224$ & $14$ & $(42,1323,6489)$\\
$F_{11}$ & $(1,1,1,1,1,1,1,4)$ & $32$ & $8$ & $(28,0,0)^*$\\
$F_{12}$ & $(1,1,1,1,1,2,3,3)$ & $672$ & $18$ & $(36,252,464)$\\
$F_{13}$ & $(1,1,1,1,2,2,2,3)$ & $1,120$ & $20$ & $(38,585,3210)$\\
$F_{14}$ & $(1,1,1,1,1,2,2,4)$ & $672$ & $18$ & $(32,66,36)$\\
$F_{15}$ & $(1,1,1,2,2,2,3,3)$ & $2,240$ & $24$ & $(33,120,302)$\\
$F_{16}$ & $(1,1,1,1,2,2,3,4)$ & $3,360$ & $30$ & $(31,62,82)$\\
$F_{17}$ & $(1,1,1,1,2,2,2,5)$ & $1,120$ & $20$ & $(25,3,0)^*$\\
$F_{18}$ & $(1,1,1,1,1,3,3,4)$ & $672$ & $18$ & $(27,0,1)^*$\\
$F_{19}$ & $(1,1,1,2,2,3,3,4)$ & $6,720$ & $36$ & $(28,22,22)$\\
$F_{20}$ & $(1,1,1,1,2,3,3,5)$ & $3,360$ & $30$ & $(25,2,1)^*$\\
$F_{21}$ & $(1,1,2,2,2,3,3,5)$ & $6,720$ & $36$ & $(24,3,1)^*$\\
$F_{22}$ & $(1,1,1,2,2,3,4,5)$ & $13,440$ & $48$ & $(24,3,1)^*$\\
\hline\hline
\end{tabular}
\end{center}
\end{table}

\begin{table}
\begin{center}
\caption{Orbits of vertices of the hypermetric polytope $HYPP_8$ originating from extreme Delaunay polytopes $2_{21}$ and $3_{21}$.
Column 2 is the order of the stabilizer;
column 5 is number of orbits of type $2_{21}$ and  $3_{21}$
that merged into single orbit $V_i$ in $HYPP_8$;
columns 6, 7 are the number of facets, containing the orbit
representative, and the number of vertices of $HYPP_8$ adjacent to it}
\label{HYP8polytope_23_21_vertices}
\begin{tabular}{||c|c|c|c|c|c|c||}
\hline
\hline
$V_{i}$ & $\vert Stab \vert$ & $\frac{\left\vert V_i\right\vert}{
10,752}$ & $\#orbits\;Sym(8)$ & Merging $23_{21}$ & Incidence & Adjacency\\
\hline
$V_{1}$ & $24$ & $20$ & $36$ & $2_{21}(8), 3_{21}(12)$ & $112$ & $848$\\
$V_{2}$ & $48$ & $10$ & $30$ & $2_{21}(7), 3_{21}(10)$ & $104$ & $799$\\
$V_{3}$ & $96$ & $5$ & $23$ & $2_{21}(5), 3_{21}(8)$ & $94$ & $701$\\
$V_{4}$ & $12$ & $40$ & $48$ & $2_{21}(11), 3_{21}(16)$ & $94$ & $758$\\
$V_{5}$ & $8$ & $60$ & $46$ & $2_{21}(10), 3_{21}(16)$ & $94$ & $804$\\
$V_{6}$ & $12$ & $40$ & $40$ & $2_{21}(8), 3_{21}(16)$ & $92$ & $979$\\
$V_{7}$ & $240$ & $2$ & $18$ & $2_{21}(2), 3_{21}(5)$ & $86$ & $926$\\
$V_{8}$ & $12$ & $40$ & $48$ & $2_{21}(11), 3_{21}(16)$ & $86$ & $709$\\
$V_{9}$ & $8$ & $60$ & $54$ & $2_{21}(12), 3_{21}(18)$ & $86$ & $728$\\
$V_{10}$ & $4$ & $120$ & $72$ & $2_{21}(16), 3_{21}(24)$ & $86$ & $774$\\
$V_{11}$ & $16$ & $30$ & $33$ & $2_{21}(6), 3_{21}(13)$ & $84$ & $1,070$\\
$V_{12}$ & $4$ & $120$ & $60$ & $2_{21}(12), 3_{21}(24)$ & $84$ & $963$\\
$V_{13}$ & $48$ & $10$ & $26$ & $2_{21}(4), 3_{21}(8)$ & $82$ & $1,023$\\
$V_{14}$ & $12$ & $40$ & $48$ & $2_{21}(7), 3_{21}(18)$ & $81$ & $1,080$\\
$V_{15}$ & $4$ & $120$ & $60$ & $2_{21}(9), 3_{21}(30)$ & $79$ & $935$\\
$V_{16}$ & $20$ & $24$ & $24$ & $2_{21}(5), 3_{21}(8)$ & $78$ & $734$\\
$V_{17}$ & $16$ & $30$ & $33$ & $2_{21}(7), 3_{21}(12)$ & $78$ & $679$\\
$V_{18}$ & $8$ & $60$ & $46$ & $2_{21}(10), 3_{21}(16)$ & $78$ & $690$\\
$V_{19}$ & $4$ & $120$ & $56$ & $2_{21}(12), 3_{21}(20)$ & $78$ & $716$\\
$V_{20}$ & $4$ & $120$ & $56$ & $2_{21}(9), 3_{21}(25)$ & $78$ & $1,050$\\
$V_{21}$ & $60$ & $8$ & $16$ & $2_{21}(2), 3_{21}(6)$ & $76$ & $1,070$\\
$V_{22}$ & $4$ & $120$ & $48$ & $2_{21}(9), 3_{21}(20)$ & $76$ & $941$\\
$V_{23}$ & $16$ & $30$ & $29$ & $2_{21}(5), 3_{21}(11)$ & $74$ & $1,032$\\
$V_{24}$ & $4$ & $120$ & $48$ & $2_{21}(8), 3_{21}(22)$ & $74$ & $920$\\
\hline\hline
\end{tabular}
\end{center}
\end{table}

\begin{table}
\begin{center}
\caption{Incidence between the orbits $V_i$ of vertices and the  orbits $F_j$ of facets of $HYPP_7$.
For each representative of orbit $V_i$, the number of incident facets is given.}
\label{HYPP7_incidence}
\begin{tabular}{||c|c|c|c|c|c|c|c||}
\hline
\hline
 & $F_{1}$ & $F_{2}$ & $F_{3}$ & $F_{4}$ & $F_{5}$ & $F_{6}$ & $F_{7}$\\
\hline
$V_{1}$ & 105 & 210 & 35 & 630 & 546 & 147 & 2100\\
$V_{2}$ & 8 & 6 & 0 & 4 & 2 & 0 & 1\\
$V_{3}$ & 11 & 6 & 1 & 2 & 1 & 0 & 0\\
$V_{4}$ & 12 & 7 & 0 & 2 & 0 & 0 & 0\\
$V_{5}$ & 15 & 5 & 1 & 0 & 0 & 0 & 0\\
$V_{6}$ & 14 & 7 & 0 & 0 & 0 & 0 & 0\\
\hline\hline
\end{tabular}
\end{center}
\end{table}

\section{Structure for all simplices}\label{Structure_AllSimplices}

The hypermetric cone $HYP_{n+1}$ describes the possible ways,in which the simplex $S$
on vertices $0$, $e_1$, \dots, $e_n$ can be embedded in a Delaunay polytope.
In particular, we saw in Section \ref{Facet_HYP8} how the facets of $HYP_{n+1}$ correspond
to the repartitioning polytopes, in which $S$ can be embedded.

In dimension $n\geq 5$ there are other simplices, necessarily of volume higher than $1$,
that can define Delaunay polytopes.
In \cite{BaranovskiDim5}, it was determined that the possible volumes of $5$-dimensional
Delaunay simplices are $1$ and $2$.
In \cite{Baranovskii95,Baranovskii99}, it was determined that the volumes of $6$-dimensional
Delaunay simplices are $1$, $2$ and $3$. Then, in \cite{RyshkovBaranovskii98}, the
facets of the corresponding Baranovskii cones $Bar_S$ were found.

In \cite{InhomogeneousPerfect}, the list of all possible $7$-dimensional Delaunay simplices
has been determined and from that we can get for each simplex its associated Baranovskii cone $Bar_S$.
There are $11$ types of such simplices and, in contrast to the lower dimensional cases, two simplices
can have the same volume and yet be inequivalent.
Key information about those Delaunay simplices are given in Table \ref{FacetOtherCones}.

\begin{table}
\begin{center}
\caption{All types of Delaunay simplices in $\mathbb{R}^7$, their volume, size of automorphism group and the number of facets of the cones $Bar_S$}
\label{FacetOtherCones}
\begin{tabular}{||c|c|c|c|c|c|c||}
\hline\hline
$S_{i}$ & vol($S_i$) & $|Aut|$ & $\#facets\;Bar_{S_i}$\\
\hline
$S_{1}$ & 1 & 40,320 & 298,592(86)\\
$S_{2}$ & 2 & 40,320 & 5,768(9)\\
$S_{3}$ & 2 & 1,440 & 6,590(62)\\
$S_{4}$ & 3 & 540 & 966(9)\\
$S_{5}$ & 3 & 1,152 & 728(9)\\
$S_{6}$ & 3 & 240 & 640(39)\\
$S_{7}$ & 4 & 1,440 & 28(3)\\
$S_{8}$ & 4 & 240 & 153(11)\\
$S_{9}$ & 4 & 144 & 131(10)\\
$S_{10}$ & 5 & 72 & 28(6)\\
$S_{11}$ & 5 & 48 & 28(8)\\
\hline\hline
\end{tabular}
\end{center}
\end{table}

\section{Hypermetrics on graphs}\label{Hypermetric_OnGraphs}

Given a graph $G$=(V,E), the notion of cut is well defined. It suffices to restrict the cut semimetric
on the edges of the graph and one obtains the {\em cut polytope of the graph} $CUTP(G)$.
The notion of $METP_n$ can also be extended to the graph setting but requires more work:
 for a cycle $C$ and an odd sized set $F$ of edges in $C$, the {\em cycle inequality} $m_{C,F}$ is
defined as
\begin{equation*}
\sum_{e\in C - F} x_e - \sum_{e\in F} x_e \leq \left\vert F\right\vert - 1.
\end{equation*}
The {\em metric polytope of the graph} $METP(G)$  is defined as the polytope defined by all cycle
inequalities $m_{C,F}$ and the {\em non-negativity inequalities} $x_e\in [0,1]$.
In fact, it is the projection of $MET_n$
on $\mathbb{R}^{|E|}$, indexed by the edges of $G$.
It is known that $METP(G)=CUTP(G)$ if and only if $G$ has no $K_5$-minor.

\begin{proposition}\label{GenerateValidInequality}
Let us take a valid inequality on $CUTP_n$ of the form
\begin{equation*}
f(x)=\sum_{1\leq i<j\leq n} a_{ij} x_{ij} \leq C.
\end{equation*}
Suppose that  we have $n$ vertices $v_1, \dots, v_n$ of $G$ with any two vertices
$v_i$, $v_j$ being  joined by a such  path $P_{ij}$  that:
\begin{itemize}
\item the edge set of all paths $P_{ij}$ are disjoint;
\item if $a_{ij}>0$, then $P_{ij}$ is reduced to an edge.
\end{itemize}
Then, the following inequality
\begin{equation*}
i_{f, G}(x)=\sum_{1\leq i<j\leq n} a_{ij} \left(\sum_{e\in P_{ij}} x_e\right) \leq C
\end{equation*}
is valid on $CUTP(G)$.
\end{proposition}
\proof Let us take a cut of $G=(V,E)$ defined by $S\subset V$.
If $S$ cuts the paths $P_{ij}$ in at most one edge, then the inequality
on $CUTP(G)$ reduces to the one on $CUTP_n$ and so, is valid.

In the general case, we will create a new cuts $S'$, which will allow
us to prove the required inequality.
If $S$ cuts $P_{ij}$ in more than one edge, then $a_{ij}\leq 0$. 
If both $v_i$ and $v_j$ are in the same part of the partition $(S, V-S)$,
then we set all vertices of $P_{ij}$ to be in the same part of the
partition $(S', V - S')$.
Otherwise, there exists an edge $e=\{w_i, w_j\}\in P_{ij}$ cut by $S$.
If $w_j\in S'$, then we set the vertices between $w_j$ and $v_j$ to belong
to $S'$ and we set the vertices from $w_i$ to $v_i$ not to belong to it.
Due to the sign condition on $a_{ij}$, one obtains
\begin{equation*}
i_{f, G}(\delta_S) \leq i_{f, G}(\delta_{S'}).
\end{equation*}
Then, from $S'$ we can obtain very simply a cut $S''$ on $K_n$ and this gives 
\begin{equation*}
i_{f, G}(\delta_{S'}) = f(\delta_{S''}) \leq C,
\end{equation*}
which proves the required inequality. \qed

When applied to the metric inequalities of $K_n$ and taking switchings, the above
proposition gives us the metric polytope $METP(G)$.
Therefore, it is temping to define the hypermetric polytope $HYPP(G)$ as the
polytope defined by the switchings of the extension of of all hypermetric inequalities
obtained from above Proposition. What is not clear is when $CUTP(G)=HYPP(G)$ and
whether there is a nice characterization of such hypermetrics. Above discussion is applied also to cones.

Any $K_n$ subgraph of $G$ will satisfy the hypothesis and the facets of $CUTP_n$
will give facets of $CUTP(G)$. Proposition \ref{GenerateValidInequality} gives
valid inequality induced by a class of homeomorphic $K_n$.
(A graph $H$ is {\em homeomorphic} to a subgraph of $G$ if $H$ can be mapped to
$G$ so that the edges of $H$ are mapped to disjoint paths in $G$.)
A graph $H$ is a {\em minor} of $G$ if $H$ can be obtained from $G$ by deleting
edges and vertices and contracting edges. An homeomorphic graph $K_n$ is a special
case of a $K_n$ minor.
In \cite{SeymourMatroidMulticommodity} it is proved that $CUTP(G)=METP(G)$ if and only
if $G$ has no $K_5$ as a minor. However, the proof appears nonconstructive and does not
seem to be able to give hypermetric inequalities, or their generalization, in a
straightforward way.

\section{Infinite hypermetrics}\label{Infinite_Hypermetric}

Another interesting question is to define infinite hypermetric cones.
One way to do that is to define $HYP_{\infty}$ by imposing that for all
$b\in \ZZ^{\infty}$ with finite {\em support} (i.e., the set $\{i: b_i\neq 0\}$) and $\sum_i b_i=1$ it holds 
\begin{equation*}
\sum_{i<j} b_i b_j d_{ij}\leq 0.
\end{equation*}

For example, the path metric of the skeleton of the infinite hyperoctahedron $K_{2, \dots, 2, \dots}$
is an infinite hypermetric, which does not embed isometrically into $l_1$.
In general, it is easy to build infinite hypermetrics; it basically suffices to use
Delaunay polytopes of infinite lattices. For example, above 
$\infty$-dimensional hyperoctahedron can be considered a Delaunay polytope in
the $\infty$-dimensional root lattice $\mathsf{D}_{\infty}$.

However, as far as we know, there is no general theory of Delaunay polytopes
in infinite dimensional lattices.

\bibliographystyle{amsplain_initials_eprint}
\bibliography{LatticeRef}

\providecommand{\bysame}{\leavevmode\hbox to3em{\hrulefill}\thinspace}
\providecommand{\MR}{\relax\ifhmode\unskip\space\fi MR }
\providecommand{\MRhref}[2]{%
  \href{http://www.ams.org/mathscinet-getitem?mr=#1}{#2}
}
\providecommand{\href}[2]{#2}
\begin{thebibliography}{10}

\bibitem{4ti2}
4ti2 team, \emph{4ti2--a software package for algebraic, geometric and
  combinatorial problems on linear spaces.}, URL: \url{http://www.4ti2.de/}.

\bibitem{OrbitFacetCutPolytope6}
D.~Avis and Mutt, \emph{All the facets of the six-point {H}amming cone},
  European J. Combin. \textbf{10} (1989), no.~4, 309--312, URL:
  \url{http://dx.doi.org/10.1016/S0195-6698(89)80002-2},
  doi:10.1016/S0195-6698(89)80002-2.

\bibitem{BaranovskiDim5}
E.~P. Baranovski{\u\i}, \emph{Volumes of {$L$}-simplexes of five-dimensional
  lattices}, Mat. Zametki \textbf{13} (1973), 771--782.

\bibitem{Baranovskii95}
E.~P. Baranovskii, \emph{About l-simplexes of $6$-dimensional lattices (in
  russian)}, Second International conference ``Algebraic, Probabilistic,
  Geometrical, Combinatorial and Functional Methods in the theory of numbers,
  1995.

\bibitem{Ba}
E.~P. Baranovski{\u\i}, \emph{The conditions for a simplex of $6$-dimensional
  lattice to be $l$-simplex (in russian)}, Ivan. Univ. \textbf{2} (1999),
  no.~3, 18--24.

\bibitem{Baranovskii99}
E.~P. Baranovskii, \emph{The conditions for a simplex of $6$-dimensional
  lattice to be $l$-simplex (in russian)}, Nauch. Trud. Ivan. (1999), no.~2,
  18--24.

\bibitem{CR}
T.~Christof and G.~Reinelt, \emph{Combinatorial optimization and small
  polytopes}, Top \textbf{4} (1996), no.~1, 1--64, With discussion, URL:
  \url{http://dx.doi.org/10.1007/BF02568602}, doi:10.1007/BF02568602.

\bibitem{RidgeGraphMetric}
A.~Deza and M.~Deza, \emph{The ridge graph of the metric polytope and some
  relatives}, Polytopes: abstract, convex and computational ({S}carborough,
  {ON}, 1993), NATO Adv. Sci. Inst. Ser. C Math. Phys. Sci., vol. 440, Kluwer
  Acad. Publ., Dordrecht, 1994, pp.~359--372.

\bibitem{OrbitVerticesMetricPolytope7}
A.~Deza, M.~Deza, and K.~Fukuda, \emph{On skeletons, diameters and volumes of
  metric polyhedra}, Combinatorics and computer science ({B}rest, 1995),
  Lecture Notes in Comput. Sci., vol. 1120, Springer, Berlin, 1996,
  pp.~112--128, URL: \url{http://dx.doi.org/10.1007/3-540-61576-8_78},
  doi:10.1007/3-540-61576-8-78.

\bibitem{OrbitPolytopeMetricPolytope8}
A.~Deza, K.~Fukuda, T.~Mizutani, and C.~Vo, \emph{On the face lattice of the
  metric polytope}, Discrete and computational geometry, Lecture Notes in
  Comput. Sci., vol. 2866, Springer, Berlin, 2003, pp.~118--128, URL:
  \url{http://dx.doi.org/10.1007/978-3-540-44400-8_12},
  doi:10.1007/978-3-540-44400-8-12.

\bibitem{DGL92}
M.~Deza, V.~P. Grishukhin, and M.~Laurent, \emph{Extreme hypermetrics and
  {$L$}-polytopes}, Sets, graphs and numbers ({B}udapest, 1991), Colloq. Math.
  Soc. J\'anos Bolyai, vol.~60, North-Holland, Amsterdam, 1992, pp.~157--209.

\bibitem{DGL95}
M.~Deza, V.~P. Grishukhin, and M.~Laurent, \emph{Hypermetrics in geometry of
  numbers}, Combinatorial optimization ({N}ew {B}runswick, {NJ}, 1992--1993),
  DIMACS Ser. Discrete Math. Theoret. Comput. Sci., vol.~20, Amer. Math. Soc.,
  Providence, RI, 1995, pp.~1--109.

\bibitem{Hyp7}
M.~Deza and M.~Dutour, \emph{The hypermetric cone on seven vertices},
  Experiment. Math. \textbf{12} (2003), no.~4, 433--440, URL:
  \url{http://projecteuclid.org/euclid.em/1087568019}.

\bibitem{DD15}
M.~Deza and M.~Dutour~Sikiri{\'c}, \emph{Enumeration of the facets of cut
  polytopes over some highly symmetric graphs}, preprint at {\tt
  arxiv:arXiv:1501.05407}, March 2013.

\bibitem{DL}
M.~Deza and M.~Laurent, \emph{Geometry of cuts and metrics}, Algorithms and
  Combinatorics, vol.~15, Springer, Heidelberg, 2010, First softcover printing
  of the 1997 original [MR1460488], URL:
  \url{http://dx.doi.org/10.1007/978-3-642-04295-9},
  doi:10.1007/978-3-642-04295-9.

\bibitem{DutourAdj}
M.~Dutour, \emph{Adjacency method for extreme {D}elaunay polytopes},
  Proceedings of ``Third Vorono\"\i \ Conference of the Number Theory and
  Spatial Tesselations'', 2009, pp.~94--101.

\bibitem{WebPageRepartitioningDim7}
M.~Dutour~Sikiri\'c, \emph{{D}elaunay polytopes classifications}, URL:
  \url{http://mathieudutour.altervista.org/PerfectCones/}.

\bibitem{InhomogeneousPerfect}
M.~Dutour~Sikiri{\'c}, \emph{{Enumeration of inhomogeneous perfect forms}}, in
  preparation.

\bibitem{SmoothnessRegularity}
M.~Dutour~Sikiri{\'c}, K.~Hulek, and A.~Sch{\"u}rmann, \emph{Smoothness and
  singularities of the perfect form and the second voronoi compactification of
  ${\mathcal a}_g$}, preprint at {\tt arxiv:arXiv:1303.5846}, March 2013.

\bibitem{EquivariantLtypeDSV}
M.~Dutour~Sikiri{\'c}, A.~Sch{\"u}rmann, and F.~Vallentin, \emph{A
  generalization of {V}oronoi's reduction theory and its application}, Duke
  Math. J. \textbf{142} (2008), no.~1, 127--164.

\bibitem{ComplexityVoronoiDSV}
M.~Dutour~Sikiri{\'c}, A.~Sch{\"u}rmann, and F.~Vallentin, \emph{Complexity and
  algorithms for computing {V}oronoi cells of lattices}, Math. Comp.
  \textbf{78} (2009), no.~267, 1713--1731, URL:
  \url{http://dx.doi.org/10.1090/S0025-5718-09-02224-8},
  doi:10.1090/S0025-5718-09-02224-8.

\bibitem{FinckePohst}
U.~Fincke and M.~Pohst, \emph{Improved methods for calculating vectors of short
  length in a lattice, including a complexity analysis}, Math. Comp.
  \textbf{44} (1985), no.~170, 463--471, URL:
  \url{http://dx.doi.org/10.2307/2007966}, doi:10.2307/2007966.

\bibitem{OrbitFacetCutPolytope7}
V.~P. Grishukhin, \emph{All facets of the cut cone {${\bf C}_n$} for {$n=7$}
  are known}, European J. Combin. \textbf{11} (1990), no.~2, 115--117, URL:
  \url{http://dx.doi.org/10.1016/S0195-6698(13)80064-9},
  doi:10.1016/S0195-6698(13)80064-9.

\bibitem{OrbitRaysMetricCone7}
V.~P. Grishukhin, \emph{Computing extreme rays of the metric cone for seven
  points}, European J. Combin. \textbf{13} (1992), no.~3, 153--165, URL:
  \url{http://dx.doi.org/10.1016/0195-6698(92)90021-Q},
  doi:10.1016/0195-6698(92)90021-Q.

\bibitem{MinkowskianSublattices}
W.~Keller, J.~Martinet, and A.~Sch{\"u}rmann, \emph{On classifying
  {M}inkowskian sublattices}, Math. Comp. \textbf{81} (2012), no.~278,
  1063--1092, With an appendix by Mathieu Dutour Sikiri{\'c}, URL:
  \url{http://dx.doi.org/10.1090/S0025-5718-2011-02528-7},
  doi:10.1090/S0025-5718-2011-02528-7.

\bibitem{BaRy}
S.~S. Ryshkov and E.~P. Baranovski{\u\i}, \emph{Repartitioning complexes in
  $n$-dimensional lattices (with full description for $n\leq 6$)}, Voronoi
  impact on modern science, Institute of Mathematics, Kyiv, 1998, pp.~115--124.

\bibitem{RyshkovBaranovskii98}
S.~S. Ryshkov and E.~P. Baranovskii, \emph{Repartitioning complexes in
  $n$-dimensional lattices (with full description for $n\leq 6$)}, Proceedings
  of conference "Voronoi impact on modern science", Book 2, 1998, pp.~115--124.

\bibitem{Ctype_original}
S.~S. Ry{\v{s}}kov and E.~P. Baranovski{\u\i}, \emph{{$C$}-types of
  {$n$}-dimensional lattices and {$5$}-dimensional primitive parallelohedra
  (with application to the theory of coverings)}, Proc. Steklov Inst. Math.
  (1978), no.~4, 140, Cover to cover translation of Trudy Mat. Inst. Steklov
  {{\bf{1}}37} (1976), Translated by R. M. Erdahl.

\bibitem{BistellarFlip}
F.~Santos, \emph{Geometric bistellar flips: the setting, the context and a
  construction}, International {C}ongress of {M}athematicians. {V}ol. {III},
  Eur. Math. Soc., Z\"urich, 2006, pp.~931--962.

\bibitem{bookschurmann}
A.~Sch{\"u}rmann, \emph{Computational geometry of positive definite quadratic
  forms}, University Lecture Series, vol.~48, American Mathematical Society,
  Providence, RI, 2009, Polyhedral reduction theories, algorithms, and
  applications.

\bibitem{SeymourMatroidMulticommodity}
P.~D. Seymour, \emph{Matroids and multicommodity flows}, European J. Combin.
  \textbf{2} (1981), no.~3, 257--290, URL:
  \url{http://dx.doi.org/10.1016/S0195-6698(81)80033-9},
  doi:10.1016/S0195-6698(81)80033-9.

\bibitem{OrbitFacetCutPolytope5}
M.~E. Tylkin (=M.~Deza), \emph{On {H}amming geometry of unitary cubes}, Soviet
  Physics. Dokl. \textbf{5} (1960), 940--943.

\bibitem{VoronoiII}
G.~Voronoi, \emph{Nouvelles applications des param\`etres continus \`a la
  th\'eorie des formes quadratiques. {D}euxi\`eme {M}\'emoire. {R}echerches sur
  les parall\'ello\`edres primitifs.}, J. Reine Angew. Math \textbf{134}
  (1908), no.~1, 198--287.

\end{thebibliography}

\end{document}